\documentstyle{article}                   %% no pdftex on arXiv
\pdfoutput=1 \pdfcompresslevel=9 \pdfadjustspacing=1
\begin{document}                          %% required for LaTeX
\def\item#1{\par\leavevmode\hangindent=\parindent\hangafter=1%
 \llap{#1\enspace}\ignorespaces}          %% required for LaTeX
\def\ifudf#1{\expandafter\ifx\csname #1\endcsname\relax}
\newif\ifpdf \ifudf{pdfoutput}\pdffalse\else\pdftrue\fi
%\ifpdf\pdfoutput=1\pdfcompresslevel=9\pdfadjustspacing=1\fi
%\ifpdf\pdfpagewidth=210mm \pdfpageheight=297mm
% \else\special{papersize=210mm,297mm}\fi
\hsize=149mm \vsize=245mm \hoffset=5mm \voffset=0mm
\parskip=1ex minus .3ex \parindent=0pt
\everydisplay={\textstyle}
\font\tbb=bbmsl10 \font\sbb=bbmsl10 scaled 700
 \font\fbb=bbmsl10 scaled 500
\newfam\bbfam \textfont\bbfam=\tbb
 \scriptfont\bbfam=\sbb \scriptscriptfont\bbfam=\fbb
\def\bb{\fam\bbfam}%
\font\tfk=eufm10 \font\sfk=eufm7 \font\ffk=eufm5
\newfam\fkfam \textfont\fkfam=\tfk
 \scriptfont\fkfam=\sfk \scriptscriptfont\fkfam=\ffk
\def\fk{\fam\fkfam}%
\def\fn[#1]{\font\TmpFnt=#1\relax\TmpFnt\ignorespaces}
\def\em{\expandafter\ifx\the\font\tensl\rm\else\sl\fi}
\def\dt{\number\day.\number\month.\number\year/\the\time}
\def\\{\hfill\break}
\newtoks\secNo \secNo={0}
\newcount\ssecNo \ssecNo=0
\def\section #1#2\par{\goodbreak\vskip 3ex\noindent
 \global\secNo={#1}%
 \global\ssecNo=0
 {\fn[cmbx10 scaled 1200]#1#2}\vglue 1ex}
\def\subsection #1.{\goodbreak\vskip 3ex\noindent
 \global\advance\ssecNo by 1
 {\fn[cmbx10]\the\secNo.\the\ssecNo.~#1}\vglue 1ex}
\def\url #1<#2>{\leavevmode
 \ifpdf\pdfstartlink
  attr {/Border [0 0 0] /C [0.7 0.7 0.7]}%
  user {/Subtype /Link /A << /S /URI /URI (#1) >>}%
  {#2}\pdfendlink
 \else #2\fi}
\def\web:#1 {\url https://#1<#1> }
\def\arx:#1 {EPrint \url https://arxiv.org/abs/#1<arXiv:#1> }
\def\href#1<#2>{\leavevmode
 \ifpdf\pdfstartlink attr {/Border [0 0 0 ]} goto name {#1}\fi
 {#2}\ifpdf\pdfendlink\fi}
\def\label@#1:#2@{\ifudf{#1}
 \expandafter\xdef\csname#1\endcsname{#2}\else
 \errmessage{label #1 already in use!}\fi}
 \label@fig.drop:Fig 1@                   %% required for LaTeX
 \label@def.kg:Geometry@
 \label@def.sg:Subgeometry@
 \label@def.sb:Symmetry breaking@
 \label@def.mg:M{\accent "7F o}bius geometry@
 \label@def.lsg:Lie sphere geometry@
 \label@def.lg:Laguerre geometry@
 \label@fig.proj:Fig 2@
 \label@fig.vect:Fig 3@
 \label@abma08:1@
 \label@bi05:2@
 \label@bl29:3@
 \label@bo94:4@
 \label@bopi96a:5@
 \label@bopi99:6@
 \label@bosu08:7@
 \label@bjl14:8@
 \label@bpp08:9@
 \label@bpp02:10@
 \label@buca07:11@
 \label@bjr10:12@
 \label@bjr12:13@
 \label@busa12:14@
 \label@bjrs14:15@
 \label@bjrs15:16@
 \label@bjr18:17@
 \label@bjpr19:18@
 \label@bcjpr20:19@
 \label@ca05:20@
 \label@ca28:21@
 \label@ca03:22@
 \label@ce18:23@
 \label@cps21:24@
 \label@ch67:25@
 \label@cu94:26@
 \label@da73:27@
 \label@da99:28@
 \label@de11:29@
 \label@djs21:30@
 \label@du22:31@
 \label@dura11:32@
 \label@ed36:33@
 \label@gede91:34@
 \label@jtz97:35@
 \label@jhp99:36@
 \label@imdg:37@
 \label@jesu07:38@
 \label@jko13:39@
 \label@jesz21:40@
 \label@jpp21:41@
 \label@hrsy12:42@
 \label@kl93:43@
 \label@ma60:44@
 \label@mc11:45@
 \label@muni99:46@
 \label@ni18:47@
 \label@pe19:48@
 \label@ppy21:49@
 \label@pi81:50@
 \label@pi85:51@
 \label@pi86:52@
 \label@ri97:53@
 \label@roya18:54@
 \label@gost92:55@
 \label@ta91:56@
 \label@ve26:57@
 \label@wo97:58@
\def\@#1:#2@{\ifpdf\pdfdest name {#1} xyz\fi {#2}}
\def\:#1:{\href#1<\ifudf{#1}??\else\csname#1\endcsname\fi>}
\newcount\refNo \refNo=0
\def\refitem#1 {\global\advance\refNo by 1
 \item{\@#1:\number\refNo@.}}
\def\R{{\bb R}} \def\RP{{\bb R}{\sl P}} \def\P{{\bb P}}
\def\rk{\mathop{\rm rk}}
%%%%%%%%%%%%%%%%%%%%%%%%%%%%%%%%%%%%%%%%%%%%%%%%%%%%%%%%%%%%%%%
\newdimen\FigWD \FigWD=0.57\hsize
\pdfximage width \FigWD {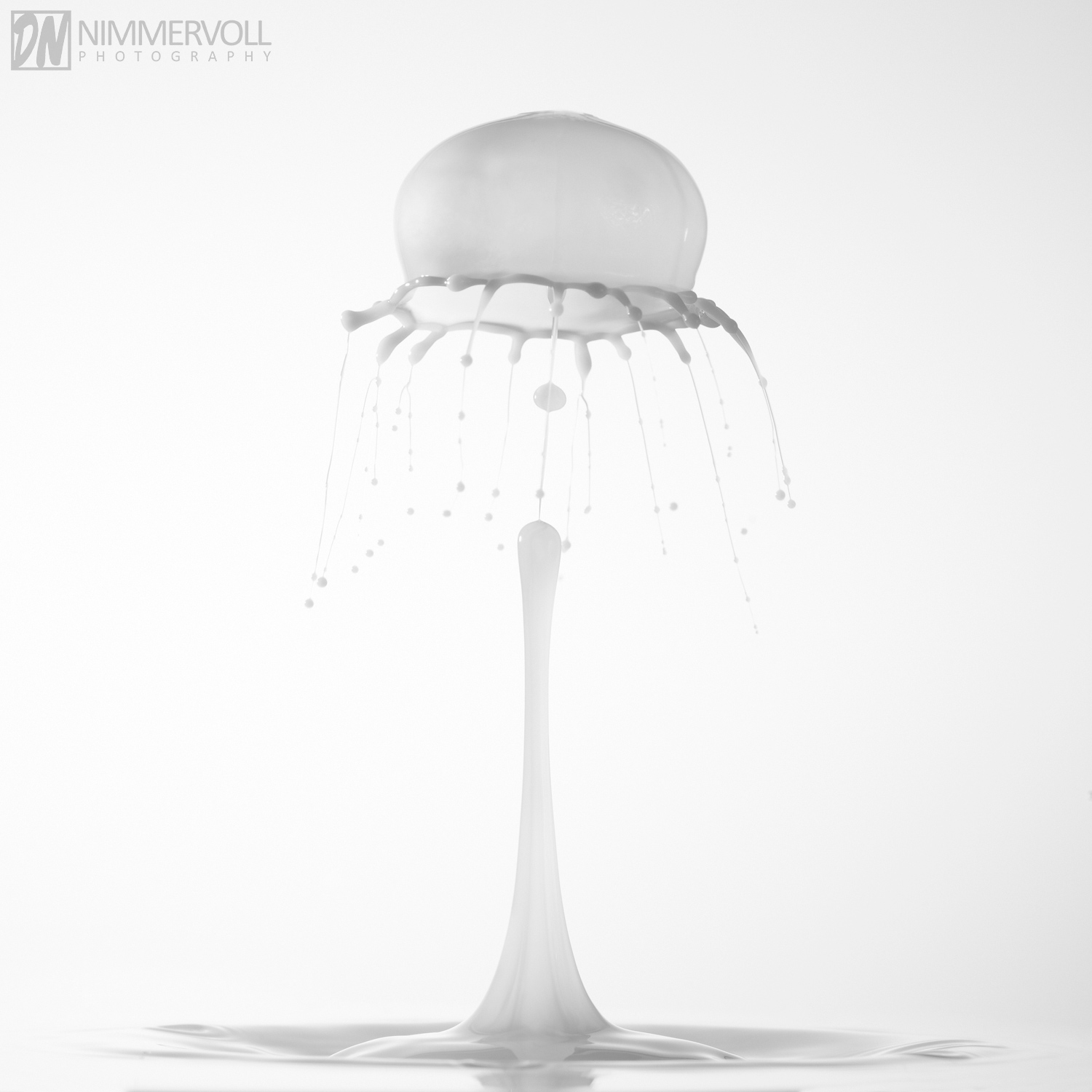}
 \xdef\DNFig{\the\pdflastximage}
\def\DoD{\vtop{\hrule height 0pt\hsize=\FigWD
 {\pdfrefximage\DNFig}\vskip 1ex
 \centerline{\vbox{\hsize=0.95\hsize\baselineskip=9pt\fn[cmr7]%
 {\fn[cmbx7]\@fig.drop:Fig 1@.}
 D Nimmervoll: Drop on Drop Splash\\
 \copyright\ 2012 Daniel Nimmervoll, \web:www.nimmervoll.org \\
 all rights reserved - re-use subject to rights-holder consent
 }}\vglue 1ex}}
%%%%%%%%%%%%%%%%%%%%%%%%%%%%%%%%%%%%%%%%%%%%%%%%%%%%%%%%%%%%%%%
\def\cascA{\vskip 2ex
\vtop{\hfil\vrule height 12pt depth 6pt width 0pt
  projective geometry of $\RP^{n+2}$ \hfil\break
 \hglue 19em
  $\swarrow$\kern 6em
  $\searrow$ \\
 \hglue 7em\vrule height 12pt depth 6pt width 0pt
  Lie sphere (contact) geometry of $S^n$ \quad\quad
  projective geometry of $\RP^{n+1}$ \\
 \hglue 12em
  $\swarrow$\kern 6em
  $\searrow$\kern 7em
  $\downarrow$ \\
 \hglue 5em\vrule height 12pt depth 6pt width 0pt
  \llap{(wider)} Laguerre geometry in $E^n$ \quad\quad
  M\"obius (conformal) geometry of $S^n$ \\
 \hglue 12em
  $\searrow$\kern 6em
  $\swarrow$\kern 4em
  $\downarrow$\kern 5em
  $\searrow$ \\
 \hglue 10.5em\vrule height 12pt depth 6pt width 0pt
  similarity geometry \quad
  hyperbolic geometry \quad
  spherical geometry
}\vglue 2ex
\centerline{\fn[cmr7]%
 {\fn[cmbx7]\@fig.proj:Fig 2@.}
 Cascade of sphere geometries: the classical viewpoint}
 \vskip 2ex}
%%%%%%%%%%%%%%%%%%%%%%%%%%%%%%%%%%%%%%%%%%%%%%%%%%%%%%%%%%%%%%%
\def\cascB{\vskip 2ex
\vtop{\hfil\vrule height 12pt depth 6pt width 0pt
  projective geometry of $\P(\R^{n+1,2})$ \hfil\break
 \hglue 19em
  $\swarrow$\kern 6em
  $\searrow$ \\
 \hglue 7em\vrule height 12pt depth 6pt width 0pt
  Lie sphere (contact) geometry of $S^n$ \quad\quad
  projective geometry of $\P(\R^{n+1,1})$ \\
 \hglue 12em
  $\swarrow$\kern 6em
  $\searrow$\kern 7em
  $\downarrow$ \\
 \hglue 5em\vrule height 12pt depth 6pt width 0pt
  \llap{(narrow)} Laguerre geometry in $E^n$ \quad\quad
  M\"obius (conformal) geometry of $S^n$ \\
 \hglue 12em
  $\searrow$\kern 6em
  $\swarrow$\kern 4em
  $\downarrow$\kern 5em
  $\searrow$ \\
 \hglue 10.5em\vrule height 12pt depth 6pt width 0pt
  Euclidean geometry \quad
  hyperbolic geometry \quad
  spherical geometry
}\vglue 2ex
\centerline{\fn[cmr7]%
 {\fn[cmbx7]\@fig.vect:Fig 3@.}
 Cascade of sphere geometries: the gauge theoretic viewpoint}
 \vskip 2ex}
%%%%%%%%%%%%%%%%%%%%%%%%%%%%%%%%%%%%%%%%%%%%%%%%%%%%%%%%%%%%%%%
\newtoks\title \newtoks\stitle \newtoks\author \newtoks\sauthor
\title={Symmetry breaking in geometry}
\stitle={Symmetry breaking}
\author={A Fuchs, U Hertrich-Jeromin, M Pember}
\sauthor={A Fuchs et al}
\ifpdf\pdfinfo{%
 /Title (\the\title) /Author (\the\author) /Date (\dt)}\fi
                                          %%  disable for LaTeX
%\headline={{\fn[cmr7]\the\sauthor\hfil\the\stitle}}
%\footline={{\fn[cmr7]\hfil-\folio-\hfil}}
\centerline{{\fn[cmbx10 scaled 1440]\the\title}}\vglue .2ex
\centerline{{\fn[cmr7]\the\author}}\vglue 3em plus 3ex
\centerline{\vtop{\hsize=.8\hsize{\bf Abstract.}\enspace
 A geometric mechanism that may,
 in analogy to similar notions in physics,
 be considered as ``symmetry breaking'' in geometry
 is described,
 and several instances of this mechanism in
 differential geometry are discussed:
 it is shown how spontaneous symmetry breaking may occur, and
 it is discussed how explicit symmetry breaking may be used
  to tackle certain geometric problems.
 A systematic study of symmetry breaking in geometry
  is proposed,
  and some preliminary thoughts on further research
  are formulated.
}}\vglue 2em
\centerline{\vtop{\hsize=.8\hsize{\bf MSC 2020.}\enspace
 {\it 53A10\/}, 53A40, 53A70, 37K35, 70H33
}}\vglue 1em
\centerline{\vtop{\hsize=.8\hsize{\bf Keywords.}\enspace
 symmetry breaking;
 conserved quantity;
 integrable reduction;
 sphere geometry;
 contact geometry;
 Moebius geometry;
 conformal geometry;
 Laguerre geometry;
 isothermic surface;
 Willmore surface;
 channel surface;
 linear Weingarten surface;
 constant curvature;
 Weierstrass representation;
 cyclic system;
 Guichard net;
 transformation;
 integrable discretization.
}}\vglue 4em
%%%%%%%%%%%%%%%%%%%%%%%%%%%%%%%%%%%%%%%%%%%%%%%%%%%%%%%%%%%%%%%

%%%%%%%%%%%%%%%%%%%%%%%%%%%%%%%%%%%%%%%%%%%%%%%%%%%%%%%%%%%%%%%
%%%%%%%%%%%%%%%%%%%%%%%%%%%%%%%%%%%%%%%%%%%%%%%%%%%%%%%%%%%%%%%
\section Introduction

In physics,
the notion of symmetry breaking is well established and has,
in some instances, attracted substantial attention.
As a consequence, the phenomenon is fairly well understood,
and even if the details of symmetry breaking may depend on
the specific physical theory it appears in, there seems to
be a general agreement about meaning and significance
of the notion.

%% image: start 2-columns
\vskip 1ex\hbox to \hsize{\vtop{\advance\hsize by -1.08\FigWD

In geometry,
on the other hand, a notion of symmetry breaking is mostly
unheard of --- despite the fact that the notion of symmetry
has become central to an understanding of geometry
(in fact, geometries),
most notably through Klein's ``Erlanger programme'' [\:kl93:],
at about the same time when symmetry considerations gained
interest in physics, cf [\:cu94:].
Thus, phenomena that could
--- in analogy to the notion in physics ---
be considered as ``symmetry breaking'' in a geometric context
were, at least to our knowledge, never investigated
in a systematic way.

The aim of this text is to give a first outline of a more
systematic approach to these phenomena,
as well as to the related geometric techniques,
in order to obtain a better understanding not only of these
phenomena but also of the interrelations between different
geometries.

 }\hfill\raise 0pt\DoD\hfil}\vskip 1ex %% end 2-columns

Thus we shall present a cabinet of curiosities of
what we consider ``symmetry breaking'' phenomena
in (differential) geometry,
trying to present a variety of ways resp contexts in which
such symmetry breaking may occur on the one hand
while, on the other hand, attempting to highlight the features
of the presented phenomena that in our view make them qualify
as ``symmetry breaking'' in order to obtain a first abstraction
for a more systematic study.

Naturally we will draw from our areas of expertise:
thus we will focus on symmetry breaking in sphere geometries,
where we or our collaborators have encountered, described
and ``collected'' many of these phenomena
--- in the case of the second author for more than three
    decades now.
Confining to sphere geometries also has the advantage
of providing a clear arrangement while offering a context
that is rich enough to display a wide variety of symmetry
breaking phenomena.

To clarify similarities and differences with symmetry breaking
in physics that we see,
we will begin with a very brief introduction of the notion
from a physics point of view,
tailored for the purpose of this text,
including an abstract description as well as several examples:
the last of those examples, illustrated in \:fig.drop:,
beautifully suggests a relationship between symmetry breaking
in physics and in geometry.
By drawing on analogies with physics in this endeavour
we not only hope to justify our terminology
but, more importantly, we hope that insight from more than
a century of relevant research in physics may help to structure
the quest into symmetry breaking in geometry as well as
to suggest the ``correct'' research questions.
However, on closer inspection,
there are also crucial differences between the notions
of ``spontaneous symmetry breaking'' in physics and geometry;
how significant those differences are may become more clear
once a definition in geometry is settled.

In the following section we will then set the scene:
we make an attempt at formulating a (working) definition
 of ``symmetry breaking'' in geometry,
 as suggested by Klein's ``Erlanger programme'' [\:kl93:];
we give a brief introduction to the sphere geometries that
 will provide the realm for our curiosity shop to be presented
 in the main section of the text; and
finally we present two examples of symmetry breaking in the
 context of elementary geometry.

With these preparations we then present a variety of instances
of what we consider as ``symmetry breaking phenomena''
in geometry
in the main section of this text, adopting three different
points of view at a more technical level.
In the first part of this section we describe phenomena
 from a classical, projective geometry viewpoint,
 closely following the ``mechanics of symmetry breaking''
 as set out in Sect~2;
most of the discussed examples concern classification problems,
 that show how ``spontaneous symmetry breaking'' may occur
 in various flavours.
In the second part we switch our viewpoint to a use
 of homogeneous coordinates that not only allows for
 a more algebraic treatment of symmetry breaking in
 the given context but, in several cases, also for
 a finer control of the symmetry breaking process.
In the final third part of the section we demonstrate how
 the additional structure provided by an integrable systems
 context provides for a very efficient description of the
 interaction between symmetry breaking and the geometry
 of the objects concerned;
here we discuss instances of explicit symmetry breaking and
 how these can be used to tackle, for example, discretization
 in differential geometry.
This last approach also suggests a close relationship between
 what we consider as ``symmetry breaking''
  in (integrable) geometry and
 integrable reductions of the underlying differential equations.

Based on this cabinet of curiosities we conclude with a first attempt
at formulating some questions or problems that may lead
to further insight into symmetry breaking in geometry.

%%%%%%%%%%%%%%%%%%%%%%%%%%%%%%%%%%%%%%%%%%%%%%%%%%%%%%%%%%%%%%%
{\it Acknowledgements.\/}
We gratefully acknowledge valuable input by our colleagues
and collaborators
V Branding,
F Burstall,
and
G Szewieczek.
Furthermore, we thank D Nimmervoll
for the permission to use his photograph \:fig.drop:
in this work free of charge.
Partial financial support for this project has been provided by
the Austrian Science Fund FWF through research project P28427
 ``Non-rigidity and Symmetry breaking'';
the third author was also supported by GNSAGA of INdAM and
the MIUR grant ``Dipartimenti die Eccellenza'' 2018--2022,
CUP:E11G18000350001, DISMA, Politecnico di Torino.

%%%%%%%%%%%%%%%%%%%%%%%%%%%%%%%%%%%%%%%%%%%%%%%%%%%%%%%%%%%%%%%
%%%%%%%%%%%%%%%%%%%%%%%%%%%%%%%%%%%%%%%%%%%%%%%%%%%%%%%%%%%%%%%
\section 1. Symmetry breaking in physics

To draw analogies as well as to detect differences
between symmetry breaking phenomena in physics and
those geometric phenomena sketched in this text later,
we give a very brief description of (our understanding of)
what symmetry breaking means in physics and illustrate
this description by some examples.
This in turn may motivate the choice of terminology
adopted in our exposition.

The phenomenon of spontaneous symmetry breaking occurs in
many different theories of classical and quantum physics.
Since quantum theories seem conceptually further away from
the differential geometry of surfaces than classical theories,
we focus on the latter.
Mathematically, the state or configuration of a classical
physical system is described by a map from a domain manifold
into a target manifold.
Typical examples of domain manifolds are the real time line
of Newtonian mechanics or Minkowski spacetime for special
relativistic theories.
The target manifold describes the degrees of freedom of the
system, which may be
 the position of a particle,
 the value of some field at a point in spacetime or
 the angle between components of a kinematic system.
The laws of the physical system under consideration are the
equations of motion, usually differential equations.
Their solutions are the physically realizable states. 

Typically, the physical laws are invariant under the
action of some Lie group on the configuration space,
for example,
 translations in time or
 the isometries of the domain or target manifold,
 when these carry a Riemannian structure.
This Lie group is then considered as a symmetry group
of the theory.

In general, the space of solutions to the equations of motion
does not consist of a single point but rather is a vector
space or a manifold, typically of twice the dimension of
the target space when the equations of motion are second
order differential equations.
Therefore, a physical state is not determined solely by the
equations of motion, but additional parameters, such as initial
conditions, are required to specify a particular state.
Even if the laws of motion have certain symmetries,
asymmetry of initial conditions leads to asymmetry of the
corresponding solution, see also [\:gede91:, \S2.1].
Consequently, the symmetries of the laws of motion are in
general not symmetries of states.
However, the symmetry group of the laws acts on the
configuration space and preserves the subspace of
phyiscally realizable states:
it maps solutions to solutions,
for further details
see [\:wo97:, Chap~3] or [\:abma08:, Chap 4].

For example, the electromagnetic field can be described
by a $2$-form field on Minkowski space time.
The homogeneous Maxwell equations confine this $2$-form
to be closed and the inhomogeneous Maxwell equations relate
its divergence to the electric $4$-current acting as a source
for the electromagnetic field.
In the source-free case, the Lorentz group, the isometry group
of Minkowski space, is a symmetry group of Maxwell's equations.
As such, it acts on the space of solutions.
On the other hand, any particular solution can be specified
by prescribing its values on a spacelike hypersurface.
The only solution that is invariant under the full Lorentz
group is the vacuum solution, where the field vanishes
everywhere.
This feature, a unique completely symmetric vacuum state,
is typical for field theories.

As another example, consider a free particle in
$3$-dimensional Euclidean space.
Its position is described by a map from the real time line
to $3$-dimensional Euclidean space;
Newton's laws that govern its motion are invariant under
the group of Euclidean motions.
The space of solutions of Newton's equations of motion is
the space of affine parametrizations of straight lines
and of constant curves, that is, points.
Any particular solution can be characterized by its position
and velocity at an arbitrary instant of time
---
thus the solution space is a $6$-dimensional space,
 the product of $3$-dimensional Euclidean space and
 a $3$-dimensional real vector space.
No solution is invariant under the full symmetry group
of Euclidean motions,
but the constant solutions have higher symmetry than the
straight line solutions:
any constant solution is invariant under the subgroup of
rotations about any axis that contains that point.
The remaining group of translations acts transitively on
these point states.
Moreover, these point states of highest symmetry are all
states of lowest energy, that is, ground states.

Although the free particle in Euclidean space is usually not
mentioned as an example of spontaneous symmetry breaking,
it definitely seems to share some properties with other
physical systems, where spontaneous symmetry breaking
is said to happen, at least at the purely classical level. 

In fact, if the symmetry is broken by means of an asymmetric
initial condition, we will speak
of {\em induced\/} or {\em explicit symmetry breaking\/}
rather than
of {\em spontaneous symmetry breaking\/}
that occurs ``spontaneously'',
i.e., without any (recognizable) external cause,
cf [\:ca03:, Sect 4.2].
This distinction will be made in several places throughout
the text that follows.

Probably the most famous example of spontaneous symmetry
breaking is the Higgs mechanism,
see for example [\:mc11:, Chap~4.4].
The electroweak theory without the Higgs field is a gauge
theory with gauge group $SU(2)\times U(1)$.
However, this theory predicts massless gauge bosons,
which are not observed in experiments.
In order to give mass to three of the four bosons,
one introduces a new field, the Higgs-field.
Its classical field equations do not admit a vacuum solution
(state of lowest energy) with the full $SU(2)\times U(1)$
symmetry: thus, no vacuum solution is invariant under the
whole symmetry group of the field equations.
Picking one such vacuum solution, describing the value of
the Higgs-field by its difference from the chosen asymmetric
vacuum solution and reformulating the field equations in terms
of this difference, hides the symmetry in the field equations.
It appears as if the $SU(2)\times U(1)$ symmetry is broken.
In particular, three of the four gauge bosons appear to
acquire mass.
It should be mentioned that although no vacuum state admits
the full gauge symmetry, there is a state of higher symmetry
that does not have the lowest energy though.

As a mechanical example for spontaneous symmetry breaking,
consider a thin elastic cylindrical rod of some length and
confine its ends to some fixed distance smaller than its
length.
If no additional forces are applied, the problem has
cylindrical symmetry, but there is no unique, symmetric
state of lowest energy.
Instead there is a $1$-parameter family of asymmetric states of
lowest energy, one of which can be transformed into any other
by a rotation around the axis through the fixed endpoints.
The most symmetric state, where the rod is squeezed to a
straight line between its endpoints still is a solution,
but it is not the solution of lowest energy and furthermore
highly unstable.

Finally, we return to the example illustrated in \:fig.drop:,
that exposes a physical phenomenon of a similar flavour
as the elastic rod example above,
 cf Edgerton's photograph [\:ed36:] and [\:gost92:, Chaps~1,3]:
the splash of a drop of fluid
(here: coloured and viscosity-altered water)
that hits a fountain created by a pair of preceding drops
falling vertically into the same fluid,
 see [\:ni18:, Sect 6.3].
Clearly, the equations of motion as well as the initial state
have an $SO(2)$ rotational symmetry about the trajectory of
the drops, as they travel towards the surface of the fluid,
as an axis.
Upon the collision of the drop and the fountain the splash
at first retains the initial rotational symmetry,
but then forms (roughly equidistant) ``tentacles'':
hence the smooth rotational symmetry is broken to yield
(approximately) a discrete rotational (dihedral) symmetry,
generated by a finite rotation about the same axis as before.

Note:
the notion of {\em spontaneous\/} symmetry breaking
in these examples relates to the loss of symmetry
in a ``physically optimal'' (least energy) solution
of the corresponding (symmetric) equations of motion
---
while this is a natural notion in physics
it is less so in geometry,
as it presupposes a variational description
of the geometric configurations under consideration.

%%%%%%%%%%%%%%%%%%%%%%%%%%%%%%%%%%%%%%%%%%%%%%%%%%%%%%%%%%%%%%%
%%%%%%%%%%%%%%%%%%%%%%%%%%%%%%%%%%%%%%%%%%%%%%%%%%%%%%%%%%%%%%%
\section 2. The mechanics of symmetry breaking in geometry

A mechanism that could be considered as ``symmetry breaking''
in geometry has been described beautifully by Klein in his
Erlanger programme [\:kl93:] of 1893,
based on a group theoretic approach to geometry [\:kl93:, \S1]:

\proclaim\@def.kg:Geometry@.
Given a manifold of geometric {\em elements\/}
and a (transformation-)group acting on it;
investigate the configurations belonging to the manifold with
respect to those properties that remain unaltered by the
transformations of the group.

This immediately induces a partial order on the geometries
on a given manifold, cf [\:kl93:, \S2]:

\proclaim\@def.sg:Subgeometry@.
A geometry $(M,H)$ is a {\em subgeometry\/} of
a geometry $(M,G)$ if the transformation group $H$ of $M$
is a subgroup of the transformation group $G$ of $M$.

Motivated by the observation, how many classical geometries
arise as subgeometries of a suitable projective geometry,
by ``adjoining'' an ``absolute configuration'' $C\subset M$,
Klein also gives an alternative characterization of
a {\em subgeometry\/} $(M,H)$ of a geometry $(M,G)$
by requiring the (transformation) subgroup $H\leq G$ to
be the stabilizer group of the absolute configuration $C$.

Thinking of a transformation group acting on an manifold
as a ``symmetry group'' of,
for example, a geometric classification problem,
a terminology of ``symmetry breaking'' suggests itself:

\proclaim\@def.sb:Symmetry breaking@.
We will speak of {\em symmetry breaking\/} when a geometric
investigation in a geometry $(M,G)$ gives rise to
a specific subgeometry $(M,H)$ or,
more specifically,
to an {\em absolute configuration\/} $C\subset M$
that fixes a subgeometry $(M,H)$ via the
stabilizer subgroup $H\leq G$ of $C$.\\
And, depending on whether $C$ is imposed or appears without
apparent cause,
we shall refer to the corresponding symmetry breaking
as {\em explicit\/} resp {\em spontaneous\/}.

Clearly, this description of our use of terminology does not
qualify as a mathematical definition,
for example, the lack of an ``apparent cause'' yields anything
but a sound mathematical criterion for symmetry breaking to
be spontaneous.
However, it provides a first abstraction for a series
of phenomena
 that we observe in geometric research and
 that seem worthy of a more systematic investigation.
Thus this description will be used as a working ``definition'',
until a better understanding of the related phenomena may lead
to a ``proper'' mathematical definition.
The aim of what follows is to discuss examples and scope
of this notion of symmetry breaking,
in order to work towards such a better understanding,
and a better definition.

Note:
our notion of {\em spontaneous\/} symmetry breaking in geometry
does {\em not\/} refer to any particular type of solution of
a geometric problem;
in fact, we will often exclude trivial or degenerate solutions
and observe ``spontaneous symmetry breaking''
as a loss of symmetry in generic solutions.

%%%%%%%%%%%%%%%%%%%%%%%%%%%%%%%%%%%%%%%%%%%%%%%%%%%%%%%%%%%%%%%
\subsection Sphere geometries.
In this text, we will focus on ``symmetry breaking''
in sphere geometries.
To this end, we will give a short introduction to the sphere
geometries that we shall discuss;
for further details the interested reader is referred to
[\:bl29:], [\:ce18:] or [\:imdg:].

We start with M\"obius (conformal) geometry as probably
the simplest and most familiar of the sphere geometries:

\proclaim\@def.mg:M\"obius geometry@.
The geometry of the hypersphere-preserving (point-)
transformations of a (conformal) sphere.

Thus the elements of M\"obius geometry are points and
hyperspheres of an $n$-sphere $S^n$ --- that needs to be
equipped with a conformal structure if one wishes to
give a differential geometric characterization of
hyperspheres as subsets of points of $S^n$.
In particular, the manifold $M$ of elements of the geometry,
that the {\em M\"obius group\/} acts on,
consists of two components,
a set of points and a set of hyperspheres,
that the M\"obius group acts on separately.
Classically, hyperspheres are {\em unoriented\/}
in M\"obius geometry;
however, in various situations it is also useful
to consider oriented hyperspheres,
e.g., when considering M\"obius geometry as a subgeometry
of Lie sphere (contact) geometry:

\proclaim\@def.lsg:Lie sphere geometry@.
The geometry of transformations of hyperspheres
(including points) of a sphere that preserve oriented contact.

Here, the manifold $M$ of elements of the geometry consists
of (oriented) hyperspheres {\em and\/} points,
as degenerate hyperspheres:
the group of Lie sphere transformations acts on hyperspheres,
and may turn a hypersphere into a point or vice versa.
An example is given by a ``parallel transformation'' that
adds a constant to the radius of every hypersphere.
The description of (hyper-)surfaces in Lie sphere geometry
relies on the notion of {\em contact elements\/},
$1$-parameter families of hyperspheres in oriented contact
(at a common touching ``point'' hypersphere).

Just as M\"obius geometry can
(by distinguishing ``point spheres'')
be obtained as a subgeometry of Lie sphere geometry,
another subgeometry is obtained by distinguishing hyperplanes
among the (oriented) hyperspheres:

\proclaim\@def.lg:Laguerre geometry@.
The geometry of transformations of hyperspheres and -planes,
respectively,
of a Euclidean space that preserve oriented contact.

Here, we think of a Euclidean space as obtained from a sphere
(via stereographic projection) by distinguishing a point
(at infinity).
In fact, fixing a ``point at infinity'' in a M\"obius geometry
yields a similarity subgeometry of M\"obius geometry rather
than a Euclidean subgeometry,
where also a unit length needs to be fixed.
Accordingly, a {\em narrow Laguerre geometry\/}
is given by the requirement that the tangential distance
of hyperspheres be preserved.

%%%%%%%%%%%%%%%%%%%%%%%%%%%%%%%%%%%%%%%%%%%%%%%%%%%%%%%%%%%%%%%
\subsection Symmetry breaking: the classical viewpoint.
Klein's Erlanger programme also provides a description of
the M\"obius and Lie sphere geometries as subgeometries
of suitable projective geometries.

Embedding the ambient sphere $S^n$ of M\"obius geometry
into a projective space $\RP^{n+1}$,
hyperspheres appear as hyperplane intersections with $S^n$.
In this way, it is readily clear that any projective
transformation $\mu$ of $\RP^{n+1}$ that fixes $S^n$,
$\mu(S^n)=S^n$,
restricts to a M\"obius transformation of $S^n$;
the converse is true, but the extension of a M\"obius
transformation to the ambient projective space requires
some more detailed work,
cf [\:kl93:, \S6] or [\:imdg:, \S1.3.9].
Thus M\"obius geometry of $S^n$ is obtained as a subgeometry
of the projective geometry of $\RP^{n+1}$ by adjoining
$S^n\subset\RP^{n+1}$ as an {\em absolute quadric\/}.

In a similar way, the Lie sphere geometry of $S^n$
is obtained as subgeometry of $\RP^{n+2}$ by adjoining
the {\em Lie quadric\/} ${\cal L}^{n+1}\subset\RP^{n+2}$
as the space of hyperspheres in $S^n$,
where ${\cal L}^{n+1}$ is a quadric that contains lines
but no higher dimensional projective subspaces:
the projective lines in ${\cal L}^{n+1}$ yield the
contact elements in Lie sphere geometry mentioned above,
cf [\:kl93:, \S7] or [\:ce18:, Sect~1.6].

\cascA

Fixing a {\em point sphere complex\/}
$S^n={\cal L}^{n+1}\cap\RP^{n+1}\subset\RP^{n+2}$
identifies M\"obius geometry as a subgeometry of Lie sphere
geometry
(cf [\:kl93:, \S7])
---
restricting projective transformations of $\RP^{n+2}$
to $\RP^{n+1}$ the orientation of hyperspheres is lost
as the polar reflection in $\RP^{n+1}\subset\RP^{n+2}$,
that reverses orientations of hyperspheres,
restricts to the identity.

We already discussed above how (the wider) Laguerre geometry
arises as a subgeometry of Lie sphere geometry
by distinguishing a ``point at infinity'',
analogous to the way in which the {\em similarity geometry\/}
of a Euclidean space arises as a subgeometry of M\"obius
geometry from the stabilizer group of a point,
that yields the centre of projection
for a stereographic projection.

{\em Hyperbolic geometry\/} can be obtained as a subgeometry
of M\"obius geometry in a similar way,
by adjoining a hypersphere (not a point) at infinity
or, equivalently, a point ``outside'' $S^n\subset\RP^{n+1}$:
by polarity, this fixes the complex of all hyperspheres that
intersect the given hypersphere orthogonally,
hence a {\em complex of hyperbolic hyperplanes\/}
in a Poincar\'e half sphere model of hyperbolic space.

{\em Spherical geometry\/}, on the other hand, is obtained
by adjoining a {\em great sphere complex\/},
that is,
a point ``inside'' $S^n\subset\RP^{n+1}$,
which can be thought of as a centre of $S^n\subset\RP^{n+1}$,
or as a centre of an antipodal reflection.

Thus we obtain the cascade of geometries of \:fig.proj:,
where arrows indicate a subgeometry relation.

%%%%%%%%%%%%%%%%%%%%%%%%%%%%%%%%%%%%%%%%%%%%%%%%%%%%%%%%%%%%%%%
\subsection Symmetry breaking: a gauge theoretic approach.
Clearly, the mechanisms to pass from a (sphere) geometry
to a subgeometry described above in terms of projective
geometry are computationally expressed in terms of homogeneous
coordinates.
However, dropping the homogeneity of coordinate vectors
allows for a more detailed geometric description
---
essentially, working with vector spaces (or bundles)
allows to identify curvature instead of just its sign.
We shall sketch these ideas below,
for details the interested reader is referred to
[\:ce18:] or [\:imdg:].

As the $n$-sphere $S^n\subset\RP^{n+1}$
and the Lie quadric ${\cal L}^{n+1}\subset\RP^{n+2}$
are described as projective light cones of $\R^{n+1,1}$
and $\R^{n+1,2}$,
respectively,
the Lie sphere and M\"obius transformation groups are
covered (given)
by the respective (projective) orthogonal groups
$O(n+1,1)$ resp $O(n+1,2)$,
cf [\:imdg:, \S1.3.14] or [\:ce18:, Sect~1.7].

Fixing a point sphere complex to pass from Lie sphere geometry
to a M\"obius subgeometry amounts to fixing the orthogonal
complement $\R^{n+1,1}\cong\{p\}^\perp\subset\R^{n+1,2}$
of a timelike line $[p]\subset\R^{n+1,2}$, to obtain the point
spheres as elements of the projective light cone in
$\{p\}^\perp$.
Here $[.]$ denotes the linear span of vectors.
Normalizing $p$, so that $(p,p)=-1$ with the inner product
$(.,.)$ on $\R^{n+1,2}$, any proper sphere admits homogeneous
coordinates $\sigma=p+s\in[p]\oplus\{p\}^\perp=\R^{n+1,2}$.
Now, since $0=(\sigma,\sigma)=-1+(s,s)$ we learn that
$s\in S^{n,1}$ and,
because reflection in $\{p\}^\perp$ reverses orientations of
spheres, $\pm s\in S^{n,1}$ represent the two orientations
of the M\"obius geometric (unoriented) hypersphere
$[s]\in\RP^{n+1}$,
cf [\:ce18:, Sect~1.5] or [\:imdg:, \S1.1.3].

Further, fixing a {\em space form vector\/} $q\in\R^{n+1,1}$
in the model space $\R^{n+1,1}\cong\{p\}^\perp$ for M\"obius
geometry the quadric
$$
  Q^n := \{y\in L^{n+1}\,|\,(y,q)=-1\},
$$
obtained as a hyperplane intersection with the light cone
$L^{n+1}\subset\R^{n+1,1}$, yields a space of constant
sectional curvature $\kappa=-(q,q)$,
see [\:imdg:, \S1.4.1].
Now, the (mean) curvature of a hypersphere $s\in S^{n,1}$
is given by $-(s,q)$, see [\:imdg:, \S1.7.9];
in particular, the hyperplanes of the quadric $Q^n$ of
constant curvature  are those spheres orthogonal to $q$ or,
otherwise said,
$S^{n,1}\cap\{q\}^\perp\subset\R^{n+1,1}$ yields the
hyperplane complex of $Q^n$.
Starting from a Lie sphere geometric setting,
the symmetry between a quadric $Q^n$ of constant curvature
and its hyperplance complex $P^n$ becomes more apparent:
$$\matrix{
  Q^n = \{y\in L^{n+2}\subset\R^{n+1,2}\,|\,
   (y,p)=0,\enspace(y,q)=-1\}, \cr
  P^n = \{y\in L^{n+2}\subset\R^{n+1,2}\,|\,
   (y,p)=-1,\enspace(y,q)=0\}. \cr
}$$
This also shows that the passage from Lie sphere geometry
to Euclidean geometry via (the narrow) Laguerre geometry
is equivalent to that via M\"obius geometry,
as the two paths only differ by the order in which the
point sphere complex and the space form vector are fixed.

Summarizing we obtain the following cascade of geometries,
arrows again indicating subgeometry relations: \cascB

%%%%%%%%%%%%%%%%%%%%%%%%%%%%%%%%%%%%%%%%%%%%%%%%%%%%%%%%%%%%%%%
\subsection Two examples from elementary geometry.
To illustrate the described concepts of symmetry breaking
we shall discuss two examples from elementary geometry:
the first example is fairly simple, but relevant to certain
integrable discretizations of surfaces,
while the second displays a fairly rich structure with
respect to symmetry breaking.

{\bf Miguel's theorem}.
In discrete differential geometry circular nets are
considered as discretizations of curvature line nets
on a surface,
cf [\:imdg:, \S8.3.16] or [\:bosu08:, Sect~3.1]:
a {\em circular net\/} is a map from a $2$-dimensional
quadrilateral cell complex into a (conformal) $n$-sphere
so that the image of each quadrilateral has a circumcircle.
As an analogue of the Ribaucour transformation of a smooth
(curvature line parametrized)
surface one says that a second circular net,
defined on the same domain,
is a {\em Ribaucour transform\/} of the first if endpoints
of corresponding edges of the two nets are concircular,
that is, corresponding edges of the two nets form
the same kind of quadrilaterals as the faces of each
net do.

Now an existence question arises:
does a given circular net $x:Z\to S^n$ admit Ribaucour
transformations for any initial point $\hat x_o$, $o\in Z$?
This existence question can be reduced to a single
quadrilateral $(ijkl)$ of the cell complex $Z$,
with edges $(ij)$, $(jk)$, $(kl)$ and $(li)$.
Namely, given $\hat x_i$ the points $\hat x_j$ and $\hat x_l$
can be chosen on the circumcircles of $x_i$, $\hat x_i$
and $x_j$ resp $x_l$;
then the remaining point $\hat x_k$ has to lie on three
circles, which is an overdetermined system.

However, it is clear that all constructed circles and points
lie on the $2$-sphere given by the face of the original net
and the initial point $\hat x_i$.
Using $x_i$ as the point at infinity to break the M\"obius
geometric symmetry, the circumcircle of the face of the
original net and the two circumcircles for edges $(ij)$
and $(il)$ form a Euclidean triangle,
with points $x_k$, $\hat x_i$ and $\hat x_l$ marked
on its sides.
By Miguel's theorem the circumcircles of each vertex of
the triangle and the marked points on the adjacent sides
intersect in a common point $\hat x_k$,
cf [\:imdg:, \S8.4.10] or [\:bosu08:, Thm~3.2].

Thus Miguel's theorem from Euclidean geometry provides
a solution of a M\"obius geometric problem,
via suitable symmetry breaking.

{\bf In- and ex-centres of a Euclidean triangle}.
Consider the following M\"obius geometric configuration
[\:jko13:]:
four points and the four circumcircles of any three
of the four points.  Thus, each pair of points lies
on two circles of an elliptic pencil; once all circles
are consistently oriented, the (unique) angle bisecting
circles in each pencil can be considered: for each point,
the three angle bisecting circles through that point belong
to an elliptic pencil, hence intersect in a second point.

Breaking symmetry by choosing one of the four points
as the point at infinity of a Euclidean plane, its three
angle bisecting circles become the angle bisectors of a
Euclidean triangle and their second intersection point
becomes the in-centre of the triangle.  Less obvious
is the fact that the other three triplets of angle
bisecting circles intersect in the ex-centres of the
triangle [\:jko13:, Lemma~3.3].

Considering the above M\"obius geometric construction
on an ellipsoid in projective space as a model for
$2$-dimensional M\"obius geometry,
it turns out that the two sets of four points are in
4-fold perspective:
there are four centres of perspective on suitable lines
joining the original points to the constructed points,
thus giving rise to a ``desmic system''.

In particular, one of the four centres of perspective
lies in the inside of the ellipsoid, hence breaking
the M\"obius geometric symmetry of the construction
to a spherical geometry with that centre of symmetry
as the centre of a round sphere:
now the four original points and the constructed in-
and ex-centres become antipodal on this round sphere.

Thus considering a simple geometric configuration
in different ambient geometries,
via symmetry breaking,
reveals a surprisingly rich structure and facilitates
a deeper understanding of the geometry of the configuration.

%%%%%%%%%%%%%%%%%%%%%%%%%%%%%%%%%%%%%%%%%%%%%%%%%%%%%%%%%%%%%%%
%%%%%%%%%%%%%%%%%%%%%%%%%%%%%%%%%%%%%%%%%%%%%%%%%%%%%%%%%%%%%%%
\section 3. Symmetry breaking in differential geometry

In differential geometry spontaneous symmetry breaking seems
to occur in a similar way as in physics,
often in the context of classification problems:
a classification problem is posed in a particular geometric
context, but the solution(s) come with additional data or
``conserved quantities'' attached,
that break the symmetry of the original geometric context
and naturally place the solution(s) in a subgeometry.

On the other hand, explicit symmetry breaking may be employed
to describe and investigate a class of geometric objects in
a rather efficient way as a subclass of a (larger) class
of objects in a ``higher'' geometry:
knowledge about the properties of objects of the ``higher''
geometry as well as about the process of symmetry breaking
does not only provide information about the objects under
investigation but also about their relations to those of
the ``higher'' geometry.

Our aim in this section is to provide evidence for both:
we shall provide short descriptions of a variety of instances
of spontaneous symmetry breaking,
focusing on how symmetry breaking occurs in each case,
and we shall give an outline of some examples that show
how explicit symmetry breaking can be useful in the study
of differential geometric (classes of) objects.

%%%%%%%%%%%%%%%%%%%%%%%%%%%%%%%%%%%%%%%%%%%%%%%%%%%%%%%%%%%%%%%
\subsection The classical differential geometry approach.

In this section we discuss various examples for our notion
of ``symmetry breaking'' that are either classical or employ
a classical point of view:
though homogeneous coordinates are used for some arguments,
the core results are independent of ambient scalings,
hence obtained from symmetry breaking in a projective setting,
cf Sect 2.2.

{\bf Thomsen's theorem\/}.
A typical instance of what we describe as ``symmetry breaking''
occurs in the classification of isothermic Willmore surfaces,
cf [\:bl29:, \S81] or [\:imdg:, Sect 3.6]:
a surface $x:\Sigma^2\to\R^3$ is {\em isothermic\/}
if it admits conformal curvature line coordinates $(u,v)$,
$$
  (dx,dx) = e^{2\varphi}(du^2+dv^2)
   \enspace{\rm and}\enspace
 -(dx,dn) = e^{2\varphi}(\kappa_1du^2+\kappa_2dv^2),
$$
with the Gauss map $n:\Sigma^2\to S^2$ of $x$,
and it is a {\em Willmore surface\/} if it is critical
for the {\em Willmore functional\/}
$$
  W(x) = \int_{\Sigma^2}(H^2-K)\,dA
  = {1\over 4}\int_{\Sigma^2}(\kappa_1-\kappa_2)^2dA.
$$
A key observation is that both properties,
isothermicity and being Willmore,
are conformally invariant hence M\"obius geometric notions.
In particular, the Willmore functional is the area functional
of the {\em central sphere congruence\/} of the surface,
that can be characterized as the (unique) enveloped sphere
congruence with second order contact with the surface
in orthogonal directions:
hence Willmore surfaces were called ``conformally minimal''
in the classical literature.

Imposing both conditions for a surface $x:\Sigma^2\to S^3$,
now considered as a surface in the conformal $3$-sphere
$S^3\subset\RP^4$,
forces its central sphere congruence $c:\Sigma^2\to\RP^4$
to take values in a fixed sphere complex $k\in\RP^4$,
that is, $\gamma\perp q=0$ for lifts $\gamma$ and $q$
of $c=[\gamma]$ and $k=[q]$, respectively.

Thus symmetry breaking occurs:
considering the surface $x$ as a surface in a metric
subgeometry determined by $q$, its central sphere congruence
becomes its tangent plane congruence.
And, as the central sphere congruence can here alternatively be
characterized as the {\em mean curvature sphere congruence\/}
of the surface
---
the enveloped congruence of spheres that have the same mean
curvature as the surface at every touching point
---
Thomsen's theorem is obtained:
{\em every isothermic Willmore surface (in M\"obius geometry)
is a minimal surface in some metric subgeometry\/}.

As every minimal surface is isothermic and a Willmore surface,
a classification result for isothermic Willmore surfaces
is obtained in this way.
It should be noted that minimality of a surface is
scale-invariant,
hence only the sign of the ambient curvature is relevant,
not its magnitude,
cf [\:bl29:, \S81].

{\bf Vessiot's theorem and Willmore channel surfaces\/}.
The classification of isothermic channel surfaces reveals
another instance of ``symmetry breaking'',
caused by a more complex absolute configuration,
cf [\:ve26:, \S34] or [\:imdg:, \S3.7.5].

A {\em channel surface\/} is the envelope of a $1$-parameter
family of spheres;
as this is a Lie sphere geometric notion, hence a sensible
notion in every subgeometry of Lie sphere geometry,
the natural and most general realm of this classification
problem is again M\"obius geometry.
The enveloped $1$-parameter family $s:\Sigma^2\to\RP^4$
of spheres, $\rk ds\leq 1$, of a surface $x:\Sigma^2\to S^3$
is necessarily a (principal) curvature sphere congruence;
hence using isothermicity to introduce conformal curvature line
parameters $(u,v)$ it only depends on one of the parameters,
say, $s=s(u)$.
It then follows
that the enveloped sphere curve is planar;
thus the polar line $a\subset\RP^4$ of this plane yields
a {\em sphere pencil\/} as an absolute configuration attached
to an isothermic channel surface:
it consists of spheres that intersect all spheres of the
enveloped sphere curve orthogonally.
Depending on the type of the sphere pencil three cases occur,
cf [\:imdg:, \S\S1.2.3 \& 1.8.5]:
{\parindent=2em
\item{$\bullet$} {\em elliptic sphere pencil\/}, where
 the spheres of the pencil intersect in a circle in $S^3$
 that can be considered as an axis of rotation of the channel
 surface $x$ in any metric subgeometry;
\item{$\bullet$} {\em parabolic sphere pencil\/}, where
 the pencil touches $S^3$ in a point $\infty\in a\cap S^3$
 at infinity of a similarity subgeometry of M\"obius geometry,
 and the surface $x$ becomes a cylinder;
\item{$\bullet$} {\em hyperbolic sphere pencil\/}, where
 the pencil hits $S^3$ in two point (spheres) $\infty^\pm$,
 one of which can be used to distinguish a similarity
 subgeometry where the pencil consists of concentric
 spheres and the surface $x$ becomes a cone,
 with the common centre of the spheres as its apex;
 considering, alternatively, both points as points at infinity
 of a hyperbolic subgeometry $x$ becomes an equivariant
 surface,
 i.e., invariant under a $1$-parameter group of isometries.
\par}
Vessiot's theorem [\:ve26:, \S34] provides the Euclidean
(similarity) versions of this classification result,
and a Lie geometric version has recently been obtained
in [\:jpp21:, Thm~4.3].

A less familiar but more direct and also useful interpretation
employs the fixed sphere pencil directly as an
absolute configuration,
which yields symmetry breaking into one of three metric
subgeometries with $4$-dimensional isometry groups,
of $S^1\times H^2$, $\R^1\times\R^2$ or $H^1\times S^2$.

A {\em Dupin cyclide\/} is an envelope of two $1$-parameter
families of (curvature) spheres.
As touching is encoded by polarity in Lie sphere geometry,
the enveloped sphere curves take values in polar $2$-planes
in $\RP^5$:
a Dupin cyclide can be encoded algebraically as a polar pair
of planes in $\RP^5$, or as an orthogonal decomposition of
$\R^{4,2}=\R^{2,1}\oplus_\perp\R^{2,1}$ of the space of
homogeneous coordinates into two Minkowski $3$-spaces.
Thus a Lie sphere geometric classification of Dupin cyclides
yields a single surface,
cf [\:pi85:, \S6].

A M\"obius geometric classification of Dupin cyclides can
then be obtained via symmetry breaking,
by choosing a point sphere complex $S^3\subset\RP^5$,
i.e., fixing $p\in\R^{4,2}$ with $(p,p)=-1$:
now, a $1$-parameter family of Dupin cyclides is obtained,
characterized by the relative position of the polar pair
of planes and the chosen point sphere complex,
cf
[\:ma60:],
[\:ce18:, Sect~3.4]
or [\:imdg:, \S1.8.8].
However, as both enveloped sphere curves are planar (circles),
similar considerations as in Vessiot's theorem lead to further
symmetry breaking (in two ways),
and hence the well known rich classification of Dupin cyclides
in Euclidean (similarity) geometry is obtained,
cf [\:du22:, Mem IV \S3] or [\:pi86:].

As any Willmore channel surface is isothermic
[\:imdg:, \S3.7.10],
the classification of Willmore channel surfaces gives again
rise to symmetry breaking in two ways:
 as isothermic Willmore surfaces and
 as isothermic channel surfaces.
It turns out that these two ways of symmetry breaking
are compatible:
any Willmore channel surface is an equivariant minimal
surface in a space form,
where tori and cylinders appear as minimal surfaces
 in hyperbolic space,
whereas Willmore surfaces of revolution appear
 as minimal surfaces in various ambient curvatures,
cf [\:ta91:], [\:muni99:] and [\:imdg:, \S3.7.15].

{\bf Cyclic Guichard nets\/}.
Rather similar instances of ``symmetry breaking'' as those
described above occur in the context
of $3$-dimensional (locally) conformally flat hypersurfaces
 in the conformal $4$-sphere $S^4$ and, closely related,
of Guichard nets in the conformal $3$-sphere $S^3$.
The geometry of a generic conformally flat hypersurface,
i.e., a hypersurface $x:\Sigma^3\to S^4$ with flat induced
 conformal structure and with three distinct principal
 curvatures,
is reflected by the intrinsic geometry of its curvature
line net $y:\Sigma^3\to\R^3$:
this forms a triply orthogonal system of surfaces $y_i=const$,
satifying the (conformally invariant)
{\em Guichard condition\/}
$$
  l_1^2 + l_2^2 = l_3^2,
   \enspace{\rm where}\enspace
  (dx,dx) = l_1^2dy_1^2 + l_2^2dy_2^2 + l_3^2dy_3^2.
$$
In particular, a generic conformally flat hypersurface
can be reconstructed uniquely from this
{\em canonical principal Guichard net\/},
cf [\:imdg:, \S2.4.6].

A {\em cyclic system\/} is a triply orthogonal system so that
the orthogonal trajectories of one family of surfaces consist
of circular arcs;
or, alternatively, two families of surfaces are (parts of)
channel surfaces.
Imposing this additional condition of cyclicity it turns out
that a partial classification causes symmetry breaking:
associated with every cyclic Guichard net in the conformal
$3$-sphere $S^3\subset\RP^4$ there is a distinguished point
$[q]\in\RP^4$ --- using this $q$ to fix a metric subgeometry,
the circular arcs form normal lines to a family of parallel
{\em linear Weingarten surfaces\/} $y_3=t$,
i.e., surfaces satisfying
$$
  0 = a_tK + 2b_tH + c_t
    = a_t\kappa_1\kappa_2 + b_t(\kappa_1+\kappa_2) + c_t.
$$
On the other hand, most linear Weingarten surfaces in
a space of constant sectional curvature can be used to
generate a cyclic Guichard net via its parallel family
of surfaces,
cf
[\:imdg:, Sect 2.6].
Note that a linear Weingarten condition is not scale invariant
but the fact that a surface be linear Weingarten is,
hence only the sign of the ambient curvature and
the {\em type\/} of the linear Weingarten surface
(the sign of the discriminant $b_t^2-a_tc_t$ of the equation)
is relevant in this context.

This instance of symmetry breaking also sets the scene for
a M\"obius geometric approach to space form analogues of
Bonnet's classical theorem on parallel constant mean and
constant Gauss curvature surfaces,
cf
 [\:jtz97:, Sect II.3] or
 [\:imdg:, Sect 2.7] and [\:jpp21:, Prop 2.8],
as well as for an approach to discrete flat fronts in
hyperbolic space,
cf [\:djs21:] and [\:jesz21:, Sect 4].

A conformally flat hypersurface in $S^4$ can be reconstructed
from its canonical principal Guichard net, as mentioned above;
consequently,
the classification of conformally flat hypersurfaces with
cyclic Guichard net naturally complements the above
classification of cyclic Guichard nets,
cf [\:jesu07:]:
it turns out that every such hypersurface comes with
a sphere pencil $a\subset\RP^5$ attached
[\:jesu07:, \S4.2].
An appropriate choice of $[q]\in a$ to break the
M\"obius geometric symmetry then yields a $1$-parameter
family of parallel hyperspheres in the $4$-dimensional
space form defined by $q$.
The orthogonal curvature leaves of the cyclic Guichard net
of the hypersurface are obtained as intersections of the
hypersurface with these hyperspheres,
where they are (extrinsically) linear Weingarten surfaces.
As in the other cases of classification by symmetry breaking
described above,
this geometric characterization obtained by symmetry breaking
yields a construction method, starting from any suitable linear
Weingarten surface in (a hyperplane of) a space form.

%%%%%%%%%%%%%%%%%%%%%%%%%%%%%%%%%%%%%%%%%%%%%%%%%%%%%%%%%%%%%%%
\subsection Using homogeneous coordinates.

A slight change of viewpoint on the classical approach
yields a refinement which admits the treatment of further
symmetry breaking phenomena that remain hidden to the
classical approach.

{\bf Surfaces of constant mean curvature\/}.
For example, Thomsen's theorem above relies on a M\"obius
geometric characterization of minimal surfaces as those
surfaces whose central sphere congruence $c=[\gamma]$ takes
values in a fixed sphere complex $k=[q]$,
the complex of planes in a quadric $Q^3$ of constant curvature
given by $q$.
A similar characterization of constant mean curvature surfaces
in space forms,
using the mean curvature of the central spheres,
relies on the use of homogeneous coordinates:
when $\gamma$ is normalized, $(\gamma,\gamma)\equiv 1$,
then $H=-(q,\gamma)$ yields the mean curvature of the central
spheres, hence of the surface.
Thus a surface $x:\Sigma^2\to S^3$ has constant mean curvature
in the quadric $Q^3\subset L^4$ of constant curvature given
by $q$
if and only if its normalized central sphere congruence
$\gamma:\Sigma\to S^{3,1}$ takes values in an affine
hyperplane parallel to $\{q\}^\perp\subset\R^{4,1}$,
that is, if $d\gamma\perp q$,
and its mean curvature is then given by $H=-(\gamma,q)$.

This characterization can be employed to show that
any surface of constant mean curvature in a space form
is isothermic and a constrained Willmore surface,
see [\:bpp02:, Sect 3.4].
The converse, that is, an analogue of Thomsen's theorem,
holds for tori [\:ri97:, Thm 35],
but fails for general isothermic constrained Willmore
surfaces, cf [\:bpp08:, Sect 4].

Another application of this M\"obius geometric characterization
of surfaces of constant mean curvature in space forms is
a characterization of discrete constant mean curvature nets
in space forms that exploits a duality between
(a suitable lift of) an isothermic surface $x:\Sigma^2\to S^3$
and its (normalized) central sphere congruence
$\gamma:\Sigma^2\to S^{3,1}$,
see [\:bjl14:, Sect 5].
This yields a description of discrete constant mean curvature
nets in space forms via (explicit) symmetry breaking.

{\bf Linear Weingarten surfaces\/}.
The above characterization of (discrete) surfaces of constant
mean curvature in space forms is a particular case of a more
general characterization of (discrete) linear Weingarten
surfaces in space forms,
see [\:bjr12:] and [\:bjr18:, \S2].
Here the key observation is that any non-tubular linear
Weingarten surface $x:\Sigma^2\to Q^3$ in a quadric of
constant curvature envelops a (possibly complex conjugate)
pair of isothermic sphere congruences
$[\gamma^\pm]:\Sigma^2\to{\cal L}^4$,
thought of as surfaces in the Lie quadric.
Hence any (non-tubular) linear Weingarten surface is
an \char 10-surface in the sense of Demoulin [\:de11:],
and an \char 10-surface is linear Weingarten in a quadric
$Q^3\subset L^5$ of constant curvature
if and only if its two enveloped isothermic sphere congruences
take values in fixed sphere complexes, $\gamma^\pm\perp q^\pm$
for some $q^\pm\in\R^{4,2}\setminus\{0\}$
with $[q^+,q^-]=[p,q]$.
Though this is a Lie-geometric characterization in the
classical sense, the recovery of the linear Weingarten
equation requires the use of homogenous coordinates.

Clearly, an \char 10-surface gives rise to an isothermic
surface in the M\"obius geometric sense if one of the
enveloped isothermic sphere congruences $[\gamma^\pm]$
consists of point spheres,
say, $\gamma^-\perp p$ for the point sphere complex
$[p]\in\RP^5$ of a M\"obius subgeometry
of Lie sphere geometry.
In this way, (discrete) constant mean curvature surfaces
arise from (discrete) linear Weingarten surfaces via
symmetry breaking,
cf [\:bjr18:, Expl 4.2].

{\bf Weierstrass representations\/}.
Besides clarifying (subclass) relations between various classes
of surfaces, more specific benefits can be drawn from those
characterizations of surface (sub-) classes via symmetry
breaking:
for example, the rich transformation theories of isothermic
and \char 10-surfaces descend to constant mean curvature resp
linear Weingarten surfaces,
reflecting the fact that the latter classes are obtained
by {\em integrable reductions\/} from the former.
This observation, in turn, allows us to put the Weierstrass type
representations for minimal, horospherical and, more generally,
Bryant or Bianchi type linear Weingarten surfaces into the
wider context of a transformation theory,
cf [\:ch67:], [\:imdg:, \S5.6.21],
[\:bjr12:, Sect 4] and [\:pe19:].
As a consequence, analogous representations for discrete
or semi-discrete counterparts from said surface classes
are obtained in a straightforward manner,
cf [\:bopi96a:, Sect 7], [\:imdg:, \S5.7.37], [\:hrsy12:],
[\:djs21:], [\:bjr18:, Sect 4], [\:roya18:] and [\:ppy21:].

%%%%%%%%%%%%%%%%%%%%%%%%%%%%%%%%%%%%%%%%%%%%%%%%%%%%%%%%%%%%%%%
\subsection A gauge theoretic approach:
 polynomial conserved quantities.
In an integrable systems context
the approach via homogeneous coordinates described in the
previous section often admits an alternative description,
as an {\em integrable reduction\/} of the underlying
differential equation(s).
Geometrically, corresponding descriptions in terms
of ``polynomial conserved quantities'',
first established in [\:buca07:],
encode our previously described symmetry breaking phenomena
in a rather concise and beautiful way,
that readily lends itself to a transparent handling of the
symmetry breaking and its associated geometric problems.

{\bf Polynomial conserved quantities\/}.
The key ingredient of this story is a loop of flat connections
associated to the geometric objects under investigation:
for special surfaces or submanifolds this may be obtained
by injection of a ``spectral parameter'' into their Weingarten
equations,
cf [\:bo94:] and [\:imdg:, Sect 3.3];
though such families of flat connections may also be described
in a more invariant way,
without the introduction of frames or other auxiliary data,
cf [\:imdg:, Sects 5.4-5] or [\:buca07:, Part~III].
The existence of these loops of flat connections is intimately
related to a corresponding transformation theory and,
in particular, non-rigidity
of the underlying geometric objects,
cf [\:da99:], [\:bi05:] and [\:ca28:]:
for example, considering the (flat) connections of the family
associated to an isothermic surface $x:\Sigma^2\to Q\subset V$
in a conformal quadric $Q$ as connections $d_t$ on the trivial
ambient (vector) space bundle $\Sigma^2\times V$,
suitable parallel sections $\hat x:\Sigma^2\to Q$ yield
the Darboux transformations of $x$,
and its conformal (Calapso) deformation is obtained
by application $x_t=g_tx$ of the gauge transformations
$g_t$ that relate the (flat) connections $d_t$ to the
trivial connection.

A {\em polynomial conserved quantity\/} in this context is
then a polynomial $p(t)=\sum_{k=0}^na_kt^k$ with coefficients
$a_k:\Sigma^2\to V$ that is parallel for the loop of
connections, $d_tp(t)=0$;
or, otherwise said, $g_tp(t)\equiv const$,
see [\:buca07:, Parts IV \& V] or [\:busa12:].

{\bf Integrable reductions for \char 10-surfaces\/}.
We already briefly discussed in the previous section
how linear Weingarten surfaces in space forms appear
as special \char 10-surfaces,
by means of two (constant) linear sphere complexes that
break the Lie sphere geometric symmetry of \char 10-surfaces
to a space form geometry.
This instance of symmetry breaking is part of a larger scheme
that we shall describe here,
mostly summarizing results from [\:bjpr19:].

A surface in Lie sphere geometry is conveniently described
in terms of two enveloped sphere congruences
$c^\pm:\Sigma^2\to{\cal L}^4$,
that can be chosen to be isothermic
(but possibly complex conjugate)
surfaces in the Lie quadric ${\cal L}^4$
in the case of an \char 10-surface.
Thus they come with their isothermic loops of flat connections
on the trivial (vector) bundle $\Sigma^2\times\R^{4,2}$
(resp its complexification),
$$
  d^\pm_t=d+t\eta^\pm,
   \enspace{\rm where}\enspace
  \eta^\pm = \gamma^\pm\wedge d\gamma^\mp
$$
for appropriate lifts $\gamma^\pm:\Sigma^2\to L^5$
of $c^\pm=[\gamma^\pm]$,
cf [\:bjpr19:, Sect 2.3].
Note that the two connections are gauge equivalent,
$d^-_t
 = \exp({t\tau})\cdot d^+_t
 = d^+_t - t\,d\tau
$
since $\tau:=\gamma^+\wedge\gamma^-$ is nilpotent
and $d\tau=t\,(\eta^+-\eta^-)$;
hence $\sigma:\Sigma^2\to\R^{4,2}$ is $d^+_t$-parallel iff
$\exp(t\tau)\,\sigma$ is $d^-_t$-parallel,
in particular,
polynomial conserved quantities may be ``exchanged'' between
the two connections at the potential cost of changing the
polynomial degree by $1$.
Similarly, we could also work with the (always real)
``middle'' connection of [\:bjpr19:, (5)],
$d^{mid}_t :
 = {1\over 2}(d^+_t+d^-_t)
 = \exp(\pm{t\over 2}\tau)\cdot d^\pm_t.
$

Obviously, if one of the isothermic sphere congruences,
say $c^-$,
takes values in a (fixed) point sphere complex $p\in\R^{4,2}$,
$c^-\perp p$,
symmetry is broken and we obtain an isothermic surface $x=c^-$
in M\"obius geometry (with central sphere congruence $c^+$).
As $\gamma^+$ and $\gamma^-$ have parallel tangent planes,
this configuration can be characterized by the fact that
$p$ is a constant conserved quantity of $c^-$,
$$
  \forall t\in\R:
  d^-_tp
    = dp
    + t\,(p,\gamma^-)\,d\gamma^+
    - t\,(p,d\gamma^+)\,\gamma^-
    = 0;
$$
for the other connections $p$ yields linear conserved
quantities,
e.g., $p(t)=p+{t\over 2}(\gamma^+\wedge\gamma^-)(p)$
for the middle connection $d^{mid}_t$;
in any case $(p(t),p(t))=-1$.
Similarly,
L(aguerre)-isothermic surfaces can be characterized as
those \char 10-surfaces,
where one of the enveloped isothermic sphere congruences,
say $c^-$,
qualifies as a (Laguerre geometric) tangent plane map, that is,
by breaking symmetry with a hyperplane complex $q\in\R^{4,2}$;
for the middle connection this yields again a linear conserved
quantity $p(t)=q+{t\over 2}(\gamma^+\wedge\gamma^-)(q)$,
which is null this time, $(p(t),p(t))=0$.

A third class of surfaces characterized in a similar
vein are Guichard surfaces (in M\"obius geometry):
following [\:ca05:] such surfaces $x:\Sigma^2\to S^3$
can be characterized by the existence of curvature
line parameters $(u,v)$ so that
(in any space form subgeometry)
$$
  \exists c\neq 0:
   c\,EG(\kappa_1-\kappa_2)^2 = E - \varepsilon G,
   \enspace{\rm with}\enspace\cases{\phantom{-}
    (dx,dx) = E\,du^2 + G\,dv^2 & and \cr
   -(dx,dn) = E\kappa_1\,du^2 + G\kappa_2\,dv^2, \cr}
$$
where the sign $\varepsilon\in\{\pm 1\}$ distinguishes two
possible types;
note how isothermic surfaces are obtained as a limiting case
$c=0$ of type $\varepsilon=+1$ Guichard surfaces.
Here the construction of the enveloped isothermic sphere
congruences (complex conjugate for type $\varepsilon=-1$)
and a linear conserved quantity for the middle connection
(with $\deg(p(t),p(t))=1$)
are somewhat more involved and the reader is referred
to [\:bjpr19:, Sect~5.2] for details.
Summarizing we obtain the following symmetry
breaking scheme in terms of linear conserved quantities,
where the (constant coefficient) polynomial
$(p(t),p(t))\in\R[t]$
distinguishes the classes:
{\em an \char 10-surface that admits a linear conserved
 quantity $p(t)$ (for its middle connection) with
{\parindent=2em
\item{$\bullet$} $(p(t),p(t))=-1$ is isothermic
 in the M\"obius subgeometry of $p(0)$,
 cf [\:bjpr19:, Thm~5.3];
\item{$\bullet$} $(p(t),p(t))=0$ is L-isothermic
 in the Laguerre subgeometry of $p(0)$,
 cf [\:bjpr19:, Thm~5.11];
\item{$\bullet$} $(p(t),p(t))=2t-1$ is a Guichard surface
 in the M\"obius subgeometry of $p(0)$,
 cf [\:bjpr19:, Thm~5.7].
}\/}

Note how the linear conserved quantity $p(t)$ and the real
polynomial $(p(t),p(t))$ encode the geometry of the symmetry
breaking in each case:
the (constant) $t^0$-coefficient of $p(t)$ provides
 the absolute configuration for symmetry breaking, and
its $t^1$-coefficient yields a particular sphere congruence
 that is attached to the geometry of the surface;
the $t^0$-coefficient of $(p(t),p(t))$ allows to read off
 the type of subgeometry, and
its linear coefficient yields information about the interplay
 between the subgeometry and this special sphere congruence.

Returning to the characterization of (non-tubular) linear
Weingarten surfaces as \char 10-surfaces whose enveloped
isothermic sphere congruences $c^\pm$ take values in
(different) fixed sphere complexes $q^\pm$,
we infer from the above discussion that the middle connection
of a linear Weingarten surface,
thought of as an \char 10-surface,
admits two linearly independent linear conserved quantities,
$p^\pm(t)=a^\pm t+q^\pm$.
In particular, as is well known, any linear Weingarten surface
in a space form is either isothermic or a Guichard surface.
In fact, using a suitable space form projection
$$
  (\xi,\nu):\Sigma^2\to Q^3\times P^3
   \enspace{\rm with}\enspace\cases{
  Q^3 = \{ y\in L^5 \,|\, (y,p) = 0, (y,q) = -1 \} \cr
  P^3 = \{ y\in L^5 \,|\, (y,p) = -1, (y,q) = 0 \} \cr}
$$
for $p,q\in[q^+,q^-]$, the linear Weingarten condition,
$0=a\,K+2b\,H+c$,
is encoded in the middle connection as well as in
its corresponding linear conserved quantities,
see [\:bjpr19:, Prop~6.4 \& Cor~6.5]:
$$
  d^{mid}_t
  = d
  + tc\,\xi\wedge d\xi
  - tb\,(\xi\wedge d\nu+\nu\wedge d\xi)
  + ta\,\nu\wedge d\nu
   \enspace{\rm and}\enspace\cases{
    p(t) = p + t\,(-b\,\xi+a\,\nu), \cr
    q(t) = q + t\,(c\,\xi-b\,\nu). \cr}
$$
Note that the characteristic polynomials,
$(p(t),p(t))=-1-2at$ and $(q(t),q(t))=-\kappa-2ct$,
depend on the choice $(p,q)$ of space form projection data
here:
polynomial conserved quantities satisfy the usual superposition
principle.
Thus in the presence of a higher dimensional bundle
of polynomial conserved quantities,
the individual characteristic polynomials have little
significance for the geometry of symmetry breaking but
the corresponding Gram matrix provides valuable information,
cf [\:bjr10:, Sect~3],
in the case at hand
$$
  G = \left(
   {(p(t),p(t))\atop(q(t),p(t))}\,
   {(p(t),q(t))\atop(q(t),q(t))}
  \right) = \left(
   {-1-2at\atop 2bt}\,
   {2bt\atop-\kappa-2ct}
  \right),
$$
where $\det G=\kappa+2(a\kappa+c)\,t+4(ac-b^2)\,t^2$ returns
the ambient space curvature as well as
the discriminant of the linear Weingarten condition,
 the sign of which determines its type resp
 whether the enveloped isothermic sphere congruences are real
 or complex conjugate.

Higher degree polynomial conserved quantities have been
investigated in [\:busa12:],
where isothermic surfaces (in M\"obius geometry) with
a quadratic conserved quantity have been shown to yield
the classical special isothermic surfaces of Darboux [\:da99:]
and Bianchi [\:bi05:]:
these surfaces generalize constant mean curvature surfaces
(which are isothermic with a linear conserved quantity)
and
exhibit astongishingly rich geometric properties.

{\bf Conserved quantities:
 spontaneous vs explicit symmetry breaking\/}.
So far, we have discussed how (sub-)classes of surfaces may
be characterized by using polynomial conserved quantities to
break the symmetry of a higher geometry into a subgeometry,
thus we have described instances of explicit rather than
of spontaneous symmetry breaking.
Obviously though, polynomial conserved quantities can be
employed to describe or detect spontaneous symmetry breaking
as well
---
however,
as the sketched theory of polynomial conserved quantities
seems to be relatively new and not widely adopted yet,
we can only reference one non-trivial example
at the time of writing [\:cps21:, Sect~6.1]:
if a Lie applicable surface has one family of spherical
curvature lines then these spherical curvature lines are
Lie sphere transforms of a constrained elastic curve in
a space form geometry.

On the other hand, polynomial conserved quantities can,
and have been, used very efficiently to describe integrable
reductions of transformation theories as well as integrable
discretizations of surface classes by (explicit) symmetry
breaking.

We first sketch how the isothermic transformations touched
upon at the beginning of this section descend to surfaces
$\xi:\Sigma^2\to Q^3$ of constant mean curvature $H$ in
a space form $Q^3$:
considering $\xi$ as an isothermic linear Weingarten surface
we have
$$
  d^{mid}_t = d + {t\over 2}\,
  (2H\,\xi\wedge d\xi + \xi\wedge d\nu + \nu\wedge d\xi)
   \enspace{\rm and}\enspace\cases{
    p(t) = p + {t\over 2}\,\xi, \cr
    q(t) = q + {t\over 2}\,(2H\,\xi+\nu); \cr}
$$
thus with $\gamma^-=\xi$ and the central sphere congruence
$\gamma^+=\gamma=\nu+H\,\xi$ we obtain an isothermic surface
in M\"obius geometry with a linear conserved quantity
(note that $\tau=\gamma^+\wedge\gamma^-=\nu\wedge\xi$):
$$
  d^-_t
  = \exp({t\tau\over 2})\cdot d^{mid}_t
  = d + t\,\xi\wedge d\gamma
   \enspace{\rm and}\enspace\cases{
    p^-(t) = \exp({t\tau\over 2})\,p(t) = p, \cr
    q^-(t) = \exp({t\tau\over 2})\,q(t) = q + t\gamma. \cr}
$$
The key to see that, resp how, the Calapso and Darboux
transformations for the isothermic surface $x=[\xi]$
descend to constant mean curvature surfaces is to
control how the isothermic loop of flat connections,
hence its (linear) conserved quantities change,
cf [\:buca07:] or [\:busa12:, Sect~3].
The Calapso transformation is readily shown to descend and,
in particular, to yield the Lawson correspondence for
constant mean curvature surfaces,
cf [\:imdg:, \S5.5.29].
The Darboux transformation, on the other hand, generically
increments the degree of a polynomial conserved quantity;
however, an orthogonality condition (at an initial point)
yields a polynomial conserved quantity of the same degree,
see [\:busa12:, Thm~3.1]:
consequently, with an appropriate choice of initial value,
the Darboux transformation of isothermic surfaces descends
to the B\"acklund transformation for surfaces of constant
mean curvature.
In fact, the presented arguments are independent of the
degree of the polynomial conserved quantity used to break
the M\"obius geometric symmetry,
hence apply to the classical special isothermic surfaces
of Darboux and Bianchi touched upon above as well,
see [\:busa12:].

The isothermic transformations directly apply to the isothermic
sphere congruences $c^\pm$ enveloped by an \char 10-surface;
to see how this yields transformations for the enveloping
\char 10-surface requires some thought:
the issue is that an \char 10-surface comes
with a pair of isothermic sphere congruences $c^\pm$,
whose Calapso or Darboux transformations need to lign up
to form a new \char 10-surface.
However, 
by the gauge equivalence $d^-_t=\exp(t\tau)\cdot d^+_t$
with $\tau=\gamma^+\wedge\gamma^-$,
hence $g^+_t=g^-_t\circ\exp(t\tau)$,
the Calapso transforms $g^-_tc^-$ and $g^+_tc^+=g^-_tc^+$
of $c^-$ and $c^+$ indeed share an enveloping
{\em Calapso transform\/} of the original \char 10-surface,
cf [\:bjpr19:, Def~2.14];
and one Darboux transform, say $\hat c^-$, is enough
to construct the Ribaucour sphere congruence $\tau\hat c^-$,
whose second envelope is the corresponding
{\em Darboux transform\/} of the \char 10-surface,
cf [\:bjpr19:, Def~2.17].
At this point, a line of arguments rather similar to
the one outlined above shows how the Calapso and Darboux
transformations for \char 10-surfaces descend to corresponding
transformations of subclasses of surfaces that are obtained
by breaking the Lie geometric symmetry using a polynomial
conserved quantity,
cf [\:bjpr19:, Sects~5 \& 6].

Note how the former instance of symmetry breaking,
where constant mean curvature surfaces are obtained
 from isothermic surfaces in M\"obius geometry, 
is obtained from the latter,
where isothermic (and Guichard surfaces) are obtained
 from \char 10-surfaces,
by a symmetry breaking process,
with a point sphere complex as a constant conserved quantity.

A core observation in integrable discretizations is that
iterated Darboux-B\"acklund transformations of a surface
generate discrete nets (via permutability properties),
as the orbit of a point,
that exhibit similar properties as their smooth counterparts:
in particular,
the smooth transformation theories are faithfully replicated
by the discrete nets obtained in this way,
see [\:bosu08:, Chap~2].
This discretization scheme can efficiently be implemented,
using polynomial conserved quantities,
if an integrable discretization of a ``higher'' surface class
is available:
this approach was first pursued in [\:bjrs14:] and [\:bjrs15:],
where integrable discretizations for surfaces of constant mean
curvature in space forms and for the classical special
isothermic surfaces are obtained from the integrable
discretization of isothermic surfaces in M\"obius geometry
of [\:bopi96a:],
generalizing previous discretizations in Euclidean space,
cf [\:jhp99:, Sect~5], [\:bopi99:, Sect~4]
or [\:bosu08:, Def~4.47];
a similar strategy can be used to obtain discretizations
for linear Weingarten surfaces,
cf [\:bjr18:],
and for Guichard, isothermic and L-isothermic surfaces
by symmetry breaking from an integrable discretization
of \char 10-surfaces,
cf [\:bcjpr20:, Sect 7.2].

%%%%%%%%%%%%%%%%%%%%%%%%%%%%%%%%%%%%%%%%%%%%%%%%%%%%%%%%%%%%%%%
%%%%%%%%%%%%%%%%%%%%%%%%%%%%%%%%%%%%%%%%%%%%%%%%%%%%%%%%%%%%%%%
\section Conclusions \& Questions

In Sects 3.1 and 3.2 we encountered several classification
results,
where an absolute configuration appeared without apparent
cause --- which ``spontaneously'' led to a symmetry breaking
in geometry in the sense described in Sect~2.
Subsequently, exploiting the integrable systems context of
the situation, we gave a refined description in terms of
polynomial conserved quantities and saw how symmetry breaking
can be used explicitly to solve certain geometric problems.

Thus it seems fairly clear {\em how\/} what we consider as
``symmetry breaking'' occurs in geometry:
we understand its ``mechanics'',
at least in a classical setting.
However, it remains unclear {\em why\/}
such (spontanteous) symmetry breaking phenomena occur:
where do these ``absolute configurations'' come from and
why are they of the given types?

Led by some apparent analogies with physics an answer might
be sought in the bifurcation theory of differential equations:
this may lead to common explanations of (spontaneous) symmetry
breaking, in geometry and physics alike
--- however, symmetry hence symmetry breaking being inherently
geometric concepts, explanations of the causes in terms
of differential equations alone seem unsatisfactory,
and are unlikely to provide an understanding of instances
such as those described in Sect~2.4.

In fact, almost all of the instances of spontaneous symmetry
breaking that we described have an integrable systems context,
hence one may ask whether or not this is coincidental or
systemic:
can all these instances of symmetry breaking be described
in terms of an integrable reduction or
of polynomial conserved quantities?
And, on the other hand, does symmetry breaking indicate
the presence of an integrable system?
For example, should we expect an integrable system behind the
configuration of in- and ex-centres of a Euclidean triangle?

As annunciated in the introduction, a better understanding
of symmetry breaking --- why spontaneous symmetry breaking
occurs and how it relates to other features of a geometric
theory --- will require further investigations and a more
systematic approach to these phenomena in geometry.
A benefit can be expected to be a better understanding
of the interrelations between different geometries,
hence of ``geometry''.

%%%%%%%%%%%%%%%%%%%%%%%%%%%%%%%%%%%%%%%%%%%%%%%%%%%%%%%%%%%%%%%
%%%%%%%%%%%%%%%%%%%%%%%%%%%%%%%%%%%%%%%%%%%%%%%%%%%%%%%%%%%%%%%
\section References

\message{References}
\bgroup\frenchspacing\parindent=2em

\refitem abma08
 R Abraham, J Marsden:
 {\it Foundations of Mechanics\/};
 AMS Chelsea Publ, Providence, Rhode Island (2008)

\refitem bi05
 L Bianchi:
 {\it  Ricerce sulle superficie isoterme e sulla deformazione
  delle quadriche\/};
 Ann Mat III 11, 93--157 (1905)

\refitem bl29
 W Blaschke:
 {\it Vorlesungen \"uber Differentialgeometrie III\/};
 Springer Grundlehren XXIX, Berlin (1929)

\refitem bo94
 A Bobenko:
 {\it Surfaces in terms of 2 by 2 matrices.
  Old and new integrable cases\/};
 in Fordy \& Wood, Harmonic maps and integrable systems
  81--127, Vieweg, Braunschweig/Wiesbaden (1994)

\refitem bopi96a
 A Bobenko, U Pinkall:
 {\it Discrete isothermic surfaces\/};
 J reine angew Math 475, 187-208 (1996)

\refitem bopi99
 A Bobenko, U Pinkall:
 {\it Discretization of surfaces and integrable systems\/};
 Oxf Lect Ser Math Appl 16, 3--58 (1999)

\refitem bosu08
 A Bobenko, Y Suris:
 {\it Discrete differential geometry. Integrable structure\/};
 Grad Stud Math 98, AMS, Providence (2008)

\refitem bjl14
 A Bobenko, U Hertrich-Jeromin, I Lukyanenko:
 {\it Discrete constant mean curvature nets in space forms:
  Steiner's formula and Christoffel duality\/};
 Discr Comput Geom 52, 612--629 (2014)

\refitem bpp08
 C Bohle, P Peters, U Pinkall:
 {\it Constrained Willmore surfaces\/};
 Calc Var Partial Differ Equ 32, 263--277 (2008)

\refitem bpp02
 F Burstall, F Pedit, U Pinkall:
 {\it Schwarzian derivatives and flows of surfaces\/};
 in Differential geometry and integrable systems (Tokyo, 2000),
 Contemp Math 308, 39--61 (2002)

\refitem buca07
 F Burstall, D Calderbank:
 {\it Conformal submanifold geometry\/};
 Manuscript (2007)
 see also
 \arx:1006.5700 (2010)

\refitem bjr10
 F Burstall, U Hertrich-Jeromin, W Rossman:
 {\it Lie geometry of flat fronts in hyperbolic space\/};
 CR 348, 661--664 (2010)

\refitem bjr12
 F Burstall, U Hertrich-Jeromin, W Rossman:
 {\it Lie geometry of linear Weingarten surfaces\/};
 CR 350, 413--416 (2012)

\refitem busa12
 F Burstall, S Santos:
 {\it Special isothermic surfaces of type $d$\/};
 J London Math Soc 85, 571--591 (2012)

\refitem bjrs14
 F Burstall, U Hertrich-Jeromin, W Rossman, S Santos:
 {\it Discrete surfaces of constant mean curvature\/};
 RIMS Kyokuroku Bessatsu 1880, 133--179 (2014)
 also \arx:0804.2707 (2008)

\refitem bjrs15
 F Burstall, U Hertrich-Jeromin, W Rossman, S Santos:
 {\it Discrete special isothermic surfaces\/};
 Geom Dedicata 174, 1--11 (2015)

\refitem bjr18
 F Burstall, U Hertrich-Jeromin, W Rossman:
 {\it Discrete linear Weingarten surfaces\/};
 Nagoya Math J 231, 55--88 (2018),
 see also \arx:1406.1293 (2014)

\refitem bjpr19
 F Burstall, U Hertrich-Jeromin, M Pember, W Rossman:
 {\it Polynomial conserved quantities of Lie applicable
  surfaces\/};
 Manuscr Math 158, 505--546 (2019)

\refitem bcjpr20
 F Burstall, J Cho, U Hertrich-Jeromin, M Pember, W Rossman:
 {\it Discrete \char 10-nets and Guichard nets\/};
 \arx:2008.01447 (2020)

\refitem ca05
 P Calapso:
 {\it Alcune superficie di Guichard e
  le relative trasformazioni\/};
 Ann Mat Pura Appl 11, 201--251 (1905)

\refitem ca28
 P Calapso:
 {\it Una nuova trasformazione delle superficie isoterme\/};
 Rend Acc Naz Lincei 8, 287--290 (1928)

\refitem ca03
 E Castellani:
 {\it On the meaning of symmetry breaking\/};
 in K Brading, E Castellani (eds) {\it Symmetries in physics\/},
 Cambridge Univ Press, Cambridge (2003) 321--334

\refitem ce18
 T Cecil:
 {\it Lie sphere geometry and Dupin hypersurfaces\/};
 CrossWorks 7, College of the Holy Cross, Worcester MA (2018)

\refitem cps21
 J Cho, M Pember, G Szewieczek:
 {\it Constrained elastic curves and
  surfaces with spherical curvature lines\/};
 \arx:2104.11058 (2021)

\refitem ch67
 E Christoffel:
 {\it Ueber einige allgemeine Eigenschaften
  der Minimumsfl\"achen\/};
 Crelle's J 67, 218--228 (1867)

\refitem cu94
 P Curie:
 {\it Sur la sym\'etrie dans les ph\'enomnes physiques.
  Sym\'etrie d'un champ \'electrique et d'un
  champ magn\'etique\/};
 J Phys, 3rd series, 3, 393--417 (1894)
 see also excerpt in
 K Brading, E Castellani (eds) {\it Symmetries in physics\/},
 Cambridge Univ Press, Cambridge (2003) 311--314

\refitem da73
 G Darboux:
 {\it Sur une classe remarquable de 
  courbes et de surfaces alg\'ebriques 
  et sur la th\'eorie des imaginaires\/};
 Gauthiers-Villars, Paris (1873)

\refitem da99
 G Darboux:
 {\it Sur les surfaces isothermiques\/};
 Ann Ec Norm Sup 16, 491--508 (1899)

\refitem de11
 A Demoulin:
  {\it Sur les surfaces $R$ et les surfaces $\Omega$\/};
  C R 153, 590--593, 705--707;
  {\it Sur les surfaces $\Omega$\/};
  C R 153, 927--929 (1911)

\refitem djs21
 J Dubois, U Hertrich-Jeromin, G Szewieczek:
 {\it Notes on flat fronts in hyperbolic space\/};
 J Geom 113:20 (2022),
 see also \arx:2106.01168 (2021)

\refitem du22
 C Dupin:
 {\it Applications de g\'eom\'etrie et de m\'echanique,
  a la marine, aux ponts et chauss\'ees, etc\/};
 Bachelier, Paris (1822)

\refitem dura11
 S Dutta, S Ray:
 {\it Bead on a rotating circular hoop:
  a simple yet feature-rich dynamical system\/};
 \arx:1112.4697 (2011)

\refitem ed36
 H Edgerton:
 {\it Milk Drop Coronet\/};
 Negative (1936(?)) MIT Museum, Cambridge MA\\
 https://\web:collections.mitmuseum.org/object2/?id=HEE-NC-36002
 (access Sep 2021)

\refitem gede91
 H Genz, R Decker:
 {\it Symmetrie und Symmetriebrechung in der Physik\/};
 Vieweg, Braunschweig (1991)

\refitem jtz97
 U Hertrich-Jeromin, E Tjaden, M Z\"urcher:
 {\it On Guichard's nets and cyclic systems\/};
 \arx:dg-ga/9704003 (1997)

\refitem jhp99
 U Hertrich-Jeromin, T Hoffmann, U Pinkall:
 {\it A discrete version of the Darboux transform
  for isothermic surfaces\/};
 Oxf Lect Ser Math Appl 16, 59--81 (1999)

\refitem imdg
 U Hertrich-Jeromin:
 {\it Introduction to M\"obius differential geometry\/};
 London Math Soc Lect Note Series 300,
  Cambridge Univ Press, Cambridge (2003)

\refitem jesu07
 U Hertrich-Jeromin, Y Suyama:
 {\it Conformally flat hypersurfaces
  with cyclic Guichard net\/};
 Int J Math 18, 301--329 (2007)

\refitem jko13
 U Hertrich-Jeromin, A King, J O'Hara:
 {\it On the M\"obius geometry of Euclidean triangles\/};
 Elem Math 68, 96--114 (2013)

\refitem jesz21
 U Hertrich-Jeromin, G Szewieczek:
 {\it Discrete cyclic systems and circle congruences\/};
 Ann Mat Pura Appl (2022),
 see also \arx:2104.13441 (2021)

\refitem jpp21
 U Hertrich-Jeromin, M Pember, D Polly:
 {\it Channel linear Weingarten surfaces in space forms\/};
 \arx:2105.00702 (2021)

\refitem hrsy12
 T Hoffmann, W Rossman, T Sasaki, M Yoshida:
 {\it Discrete flat surfaces and linear Weingarten surfaces in
  hyperbolic $3$-space\/};
 Trans AMS 364, 5605--5644 (2012)

\refitem kl93
 F Klein:
 {\it Vergleichende Betrachtungen \"uber
  neuere geometrische Forschungen\/};
 Math Ann 43, 63--100 (1893),
 see also \arx:0807.3161 (2008)

\refitem ma60
 A Mannheim:
 {\it Application de la transformation par rayons
  vecteurs r\'eciproques \`a l'\'etude de la surface
  enveloppe d'une sph\`ere tangente \`a trois sph\`eres
  donn\'ees\/};
 Nouv Ann Math I 19, 67--79 (1860)

\refitem mc11
 G McCabe:
 {\it The Structure and Interpretation of
  the Standard Model\/};
 %% Philosophy and Foundations of Physics,
 Elsevier, Amsterdam  (2011)

\refitem muni99
 E Musso, L Nicolodi:
 {\it Willmore canal surfaces in Euclidean space\/};
 Rend Ist Mat Univ Trieste 31, 177--202 (1999)

\refitem ni18
 D Nimmervoll:
 {\it High-Speed Photography\/};
 English Edition E-Book, self-published (2018);\\
 transl of {\it Highspeedfotografie\/},
  mitp Verlag, Frechen (2018)

\refitem pe19
 M Pember:
 {\it Weierstrass-type representations\/};
 Geom Dedicata 204, 299--309 (2020)

\refitem ppy21
 M Pember, D Polly, M Yasumoto:
 {\it Discrete Weierstrass-type representations\/};\\
 \arx:2105.06774 (2021)

\refitem pi81
 U Pinkall:
 {\it Dupin'sche Hyperfl\"achen\/};
 PhD thesis, Freiburg i Br (1981)

\refitem pi85
 U Pinkall:
 {\it Dupin cyclides\/};
 Math Ann 270, 427--440 (1985)

\refitem pi86
 U Pinkall:
 {\it Dupinsche Zykliden\/};
 in Fischer: Mathematische Modelle;
 Vieweg, Braunschweig (1986)

\refitem ri97
 J Richter:
 {\it Conformal maps of a Riemannian surface
  into the space of quaternions\/};
 PhD thesis, TU Berlin (1997)

\refitem roya18
 W Rossman, M Yasumoto:
 {\it Semi-discrete linear Weingarten surfaces with
  Weierstrass-type representations and their singularities\/};
 Osaka J Math 57, 169--185 (2020)

\refitem gost92
 I Stewart, M Golubitsky:
 {\it Fearful symmetry\/};
 Dover, New York (1992)

\refitem ta91
 P Tapernoux (Supervisor K Voss):
 {\it Willmore-Kanalfl\"achen\/};
 Diploma thesis, ETH Z\"urich (1991)

\refitem ve26
 M Vessiot:
 {\it Contribution \`a la g\'eom\'etrie conforme,
  th\'eorie des surfaces\/};
 Bull Soc Math France
  54, 139--179 (1926) and 55, 39--79 (1927)

\refitem wo97
 N Woodhouse:
 {\it Geometric Quantization\/};
 %% Oxford Math Monographs,
 Clarendon Press, Oxford (1997)

\egroup
%%%%%%%%%%%%%%%%%%%%%%%%%%%%%%%%%%%%%%%%%%%%%%%%%%%%%%%%%%%%%%%
\vskip3em\vfill
\bgroup\fn[cmr7]\baselineskip=8pt
\hbox to \hsize{\hfil
\vtop{\hsize=.26\hsize
 A Fuchs\\
 Rechenraum GmbH\\
 Gartengasse 21/3\\
 1050 Vienna (Austria)\\
  andreas.fuchs@rechenraum.com
 }\hfil
\vtop{\hsize=.29\hsize
 U Hertrich-Jeromin\\
% Institute of Discrete Mathematics\\ \hglue 1em and Geometry\\
 Vienna University of Technology\\
 Wiedner Hauptstra\ss{}e 8--10/104\\
 1040 Vienna (Austria)\\
  udo.hertrich-jeromin@tuwien.ac.at
 }\hfil
\vtop{\hsize=.24\hsize
 M Pember\\
% % Dipartimento di Matematica\\
% Politecnico di Torino\\
% Corso Duca degli Abruzzi 24\\
% 10129 Torino (Italy)\\
%  mason.pember@polito.it
 Florida State University\\ \hglue 1em London Study Centre\\
 99 Great Russell Street\\
 London, WX1B 3LH (UK)\\
  mason.j.w.pember@bath.edu
 }\hfil}
\egroup
%%%%%%%%%%%%%%%%%%%%%%%%%%%%%%%%%%%%%%%%%%%%%%%%%%%%%%%%%%%%%%%
\end{document}                            %% required for LaTeX
%\bye